\documentclass[draft,11pt,oneside,intlimits]{amsart}
\usepackage[latin2]{inputenc}
\usepackage{amsmath}
\usepackage{amsthm}
\usepackage{amssymb}
\usepackage[mathscr]{eucal}
\usepackage{enumerate}

\def\CC{\mathbb{C}}
\def\KK{\mathbb{K}}

\def\NN{\mathbb{N}}
\def\RR{\mathbb{R}}
\def\a{\mathcal{A}}
\def\r{\mathcal{R}}

\newcommand*{\ws}{weak$^*$}

 \DeclareMathOperator*{\diam}{diam}
\DeclareMathOperator*{\dd}{{d\hskip-1pt}}
\def\nm{\|\kern-1pt|}
 
\newcommand*{\bC}{{\mathbb C }} 
 \newcommand*{\fM}{{\mathfrak{M}} }
\newcommand*{\fMe}{\ensuremath{ {\mathfrak{M}}_1 }}

\newcommand*{\ve}{\varepsilon}

\def\sstar{\#}

\def\cconv{\mathop{\mathrm{\overline{conv}}}}

\theoremstyle{plain}
\newtheorem{proposition}{Proposition}
\newtheorem{theorem}[proposition]{Theorem}
\newtheorem{lemma}[proposition]{Lemma}
\newtheorem{corollary}[proposition]{Corollary}

\theoremstyle{definition}
\newtheorem{definition}[proposition]{Definition}

\theoremstyle{definition}

\newtheorem{remark}[proposition]{Remark}

\DeclareMathOperator*{\supp}{supp} 

\begin{document}

\title[Rendezvous numbers in normed spaces]{Rendezvous numbers in normed spaces}

\author[B.~Farkas]{B\'alint Farkas}
\address{Technische Universit\"at Darmstadt, Fachbereich Mathematik, AG4 \newline Schlo\ss gartenstra\ss e 7, D-64289 Darmstadt, Germany }
\email{farkas@mathematik.tu-darmstadt.de}

\author[Sz.~Gy.~R{\'e}v{\'e}sz]{Szil\'ard Gy\"orgy R\'ev\'esz}
\address{Alfr\'ed R\'enyi Institute\newline
Hungarian Academy of Sciences\newline Re\'altanoda u.~13--15\newline H-1053, Budapest, Hungary}
\email{revesz@renyi.hu}

\date{}

\thanks{
The present publication was supported by the Hungarian-French Scientific and Technological
Governmental Cooperation, project \# F-10/04 and the Hungarian-Spanish Scientific and Technological
Governmental Cooperation, project \# E-38/04.}

\begin{abstract}
In previous papers, we used abstract potential theory, as developed by Fuglede and Ohtsuka, to
a systematic treatment of rendezvous numbers. We considered Chebyshev constants and energies as
two variable set functions, and introduced a modified notion of rendezvous intervals which
proved to be rather nicely behaved even for only lower semicontinuous kernels or for not
necessarily compact metric spaces.

Here we study the rendezvous and average numbers of possibly infinite dimensional normed
spaces. It turns out that very general existence and uniqueness results hold for the modified
rendezvous numbers in \emph{all} Banach spaces. We also observe the connections of these
magical numbers to Chebyshev constants, Chebyshev radius and entropy. Applying the developed
notions with the available methods we calculate the rendezvous numbers or rendezvous intervals
of certain concrete Banach spaces. In particular, a satisfactory description of the case of
$L_p$ spaces is obtained for all $p>0$.
\end{abstract}



\maketitle

\section{Introduction}\label{sec:Intro}

It was shown by O.~Gross \cite{Gross} that for a compact, connected metric space $(X,d)$ there
exists a unique number $r:=r(X)\in\RR$ such that for each finite point system $P=\{x_1,x_2,\dots,
x_n\}\subset X$, $n\in \NN$ one finds a point $x\in X$ with the average distance to $P$ being
exactly $r$, that is
\begin{equation}\label{eq:rand}
\frac{1}{n}\sum_{i=1}^n d(x,x_i)=r~.
\end{equation}
This number $r(X)$ is called the \emph{rendezvous number} of the space $X$. Using this strict
definition, that is requiring the very existence of a point $x$ with exact equality under
\eqref{eq:rand}, it is well-known that both compactness and connectedness are crucial
assumptions. However, one can relax on the requirements considering so-called \emph{weak
rendezvous numbers}, meaning that there exist two points $y,z\in X$ with their average
distances to the points $x_j$ being less or equal than $r$ and greater or equal than $r$,
respectively, \cite{Th}. Clearly, for connected spaces this is equivalent to the existence of a
strong rendezvous number. Hence it is not surprising that one can prove the existence and
uniqueness of such weak rendezvous numbers under the hypothesis of compactness, see e.g.,
\cite{Th}. However, dropping the compactness condition one can not expect uniqueness as, for
example, the case of $C(K)$ spaces shows \cite{GVV, Lin}. Furthermore, sticking to connected
spaces but relaxing on compactness is also insufficient to prove the existence of rendezvous
numbers. For example, the unit sphere of $\ell_p$ spaces have no rendezvous number (unless
$p=2,\:+\infty$) \cite{hinrichs}, \cite{hinrichs/wenzel}, \cite{Lin}, \cite{W2}.

In \cite{FRpotenc} we employed a systematic potential theoretic approach to rendezvous numbers
and introduced a modified definition of these numbers, considering also closure of the
occurring average distance sets in the construction. In the classical case of compact sets and
continuous kernels (e.g., distances on compact spaces), closure is superfluous as a continuous
image of a compact set is also compact, hence closed. However, in the more general case of
non-compact sets, like unit spheres of infinite dimensional Banach spaces, and also in case of
more general, only lower semicontinuous kernels, this approach provides its yield. In
particular, with the new definition we have found very general existence results, far beyond
the setting of metric spaces. It turned out that abstract potential theory on locally compact
spaces with a lower semicontinuous kernel is an appropriate framework for such investigations.
In \cite{FRpotenc} it was also indicated that the local compactness assumption on the space $X$
is not necessary and the results go through to metric, but not necessarily locally compact
spaces as well.

We analysed further consequences of this approach in the context of metric spaces in \cite{FRmetr},
extending and explaining a good deal of previous knowledge. In the present paper we continue the
study of rendezvous- and average numbers in normed spaces.

Let us fix some notation and introduce the necessary notions. In the abstract
potential theory developed by Fuglede \cite{Fug} and Ohtsuka \cite{Oh} the usual
assumptions are the following. $X$ is a locally compact, Hausdorff space and
$k:X\times X\rightarrow \RR_+\cup\{+\infty\}$ is a lower semicontinuous,
symmetric, positive \emph{kernel} function. Nevertheless, we will consider
possibly infinite dimensional Banach spaces, so the local compactness
assumption needs to be relaxed. We will accomplish this task on the cost of
allowing special kernel functions only, such as $k(x,y):=\|x-y\|$, which is
just the usual kernel appearing in connection with rendezvous numbers.

\label{sec:pot}
\begin{definition} \label{def:ch} For arbitrary $H,L \subset X$ the \emph{general
$n^{\text{th}}$ Chebyshev constant} of $L$ with respect to $H$ is defined as
$$ M_n(H,L):=
\sup_{w_1,\ldots,w_n \in H} \inf_{x\in L}
    {\frac{1}{n}}\Bigg( \sum_{k=1}^{n} k(x,w_k) \Bigg)~.
$$
 and the \emph{$n^{\text{th}}$ general dual Chebyshev constant} of $L$ relative to $H$ is
$$
    \overline{M}_n(H,L):= \inf_{w_1,\ldots,w_n \in H} \sup_{x\in L}
    {\frac{1}{n}}\Bigg(
   \sum_{j=1}^{n} k(x,w_j)
   \Bigg)~.
$$
\end{definition}
The first part of the definition is due to Ohtsuka \cite{Ohtsuka}. By standard considerations,
just as in the case of classical Chebyshev constants, one sees that $M_n(H,L)$ and
$\overline{M}_n(H,L)$ converge to some $M(H,L)$, $\overline{M}(H,L)\in[0,+\infty]$ (see, e.g.,
\cite{FRpotenc}, \cite{Fek} or \cite{Ohtsuka}). Furthermore
$$
    \sup_{n\in\NN}M_n(H,L)=\lim_{n\rightarrow\infty}M_n(H,L)\qquad\mbox{and}\quad\inf_{n\in\NN}\overline{M}_n(H,L)=\lim_{n\rightarrow\infty}\overline{M}_n(H,L)~.
$$
The limits $M(H,L)$, $\overline{M}(H,L)$ above are called the \emph{Chebyshev constant} and the
\emph{dual Chebyshev constant} of $L$ with respect to $H$.

If $X$ is a Hausdorff topological space, let us denote by $\fM(X)$ the set of positive, regular
Borel measures on $X$ and by $\fMe(X)$ the subset of probability measures. The notation
$\fMe^{\sstar}(X)$ is used for probability measures with finite support. Given a set
$H\subseteq X$, for the subfamily of measures concentrated on $H$ (or supported on $H$, in case
$H$ is closed, cf.~\cite[pp.~144--146]{Fug}) we use the analogous notations $\fM(H)$, $\fMe(H)$
and $\fMe^{\sstar}(H)$, respectively. The \emph{potential} of a measure $\mu\in \fM(X)$ is
$$
    U^\mu(x) := \int_X  k(x,y)\dd\mu(y) ~.
$$
In the classical potential theoretic literature various notions of
\emph{energies} appear.
Already Fuglede \cite{Fuglede} and Ohtsuka \cite{Ohtsuka} introduced the following two-variate
versions of energies (see also \cite{FRpotenc}).
\begin{definition}\label{def:qq}
Let $H,L\subset X$ be fixed, and $\mu\in\fMe(X)$ be arbitrary. First put
\begin{equation}\label{d:QmuH}
Q(\mu,H):=\sup_{x\in H} U^{\mu}(x)~,\qquad\text{and also}\qquad \underline{Q}(\mu,H):=\inf_{x\in H}
U^{\mu}(x)~.
\end{equation}
Then the \emph{quasi-uniform energy} and \emph{dual quasi-uniform energy} of $L$ with respect to
$H$ are
\begin{equation}\label{qstarcapacity}
q(H,L):= \inf_{\mu\in\fMe(H)} Q(\mu,L) \qquad \text{and}\qquad
\underline{q}(H,L):=\sup_{\nu\in\fMe(H)} \underline{Q}(\nu,L)~.
\end{equation}
\end{definition}
We use the notation $M(H):=M(H,H)$, $\overline{M}(H):=\overline{M}(H,H)$,
$\underline{q}(H):=\underline{q}(H,H)$ and $q(H):=q(H,H)$ for the diagonal (classical) cases of the quantities
given in Definitions \ref{def:ch} and \ref{def:qq}.

\begin{remark}
It is not surprising that in general the quantities $M(H), q(H)$
etc. do not posses any monotonicity properties as functions of the
set $H$. The worst consequence of this is the lack of good ``inner
regularity'' properties, for example $q(H)=\inf_{K\Subset H}q(K)$
fails to hold (we use the abbreviation $K\Subset H$ to express
that $K$ is a compact subset of $H$). However, fixing one variable
the functions $M(H,L)$, $\overline{M}(L,H)$, $\underline{q}(H,L)$
and $q(L,H)$ are increasing with respect to $H$ and decreasing
with respect to $L$, and the above mentioned problem disappears.
This particularly explains the relevance and importance of the
above two-variable definitions to our subject, see also
\cite{FRpotenc}.
\end{remark}

\begin{definition}\label{def:randi}For arbitrary subsets $H, L \subset X$ the $n^{\text{th}}$
\emph{(extended) rendezvous set} of $L$ with respect to $H$ is
\begin{align}\label{nthrandi}
R_n(H,L)&:=\bigcap_{w_1,\dots,w_n\in H} \cconv \Bigl\{p_n(x):=\frac 1n \sum_{j=1}^n k(x,w_j)~:~
x\in L \Bigr\},\\
\notag R_n(H)&:=R_n(H,H)~.
\end{align}
Correspondingly, one defines
\begin{alignat}{3}
\label{randidef} R(H,L)&:=\bigcap_{n=1}^\infty R_n(H,L)~,&\quad&& R(H)&:=R(H,H)~.
\intertext{Similarly, one defines the \emph{(extended) average set} of $L$ with respect to $H$ as}
\label{avridef} A(H,L)&:=\bigcap_{\mu\in\fMe(H)} \cconv \Bigl\{U^{\mu}(x)~~:~~ x\in L \Bigr\}
~,&\quad&& A(H)&:=A(H,H)~.
\end{alignat}
\end{definition}

\begin{remark}\label{AmH} Denoting the interval
\begin{equation}\label{Amuint}
A(\mu,L):=[\underline{Q}(\mu,L), Q(\mu,L)]= \cconv\{U^{\mu}(x)~:~x\in L\}~,
\end{equation}
we see that $R_n(H,L)$, $R(H,L)$ and $A(H,L)$ are all of the form $\bigcap_{\mu} A(\mu,H)$, with
$\mu$ ranging over all averages of $n$ Dirac measures at points of $H$, over $\fMe^{\sstar}(H)$ and
over all of $\fMe(H)$, respectively.
\end{remark}

\begin{remark}\label{sdc}
Let us explain how the above notions relate to the usual definitions of rendezvous numbers or average numbers. Suppose
that $(X,k)$ is a metric space and that the set $L$ is compact. Then there is no need  for the
closure in the above definitions, since in this case the potential $U^{\mu}$ is continuous, so the
set $A(\mu,L)$, being the continuous image of the compact set $L$,  is compact. This means that a
number $r\in\RR_{+}$ belongs to $R(H,L)$ if and only if for any finite system of (not necessarily
distinct) points $x_1,\dots,x_n \in H$ one finds points $y,z\in L$ satisfying
\begin{equation}\label{nthweakrendezvous}
\frac 1n \sum_{j=1}^n k(y,x_j)\le r \qquad\text{and}\qquad \frac 1n \sum_{j=1}^n k(z,x_j) \ge r~.
\end{equation}
This is the usual definition of weak rendezvous numbers in metric spaces (see \cite{Th}). In the
next step, we can assume that $L$ is connected. In this case, \eqref{nthweakrendezvous} is further
equivalent to the existence of a ``rendezvous point'' $x\in L$ with
\begin{equation}\label{nthstrongrendezvous}
\frac 1n \sum_{j=1}^n k(x,x_j) = r~.
\end{equation}
Of course in the above reasoning an arbitrary probability measure $\mu$ can replace the average of
Dirac measures. To sum up, for compact and connected sets $L$ of metric spaces, $R(L)$ (and $A(L)$)
is a single point, and it is the classical rendezvous (or average) number of $L$ (results of Gross
\cite{Gross}, Elton and Stadje \cite{Stadje}). For further discussion and examples see
\cite{FRpotenc}.
\end{remark}

From the above definitions it is easy to identify the lower and upper endpoints of the rendezvous
and the average intervals (see \cite{FRpotenc}).
\begin{proposition}\label{prop:reni}For arbitrary subsets $H, L \subset X$
we have
\begin{alignat}{3}\label{randneqchebyn}
R_n(H,L)&= [M_n(H,L),\overline{M}_n(H,L)]~,\quad\quad& R_n(H)&= [M_n(H),\overline{M}_n(H)]~,\\
R(H,L)&= [M(H,L), \overline{M}(H,L)]~, \quad\quad &R(H)&= [M(H), \overline{M}(H)]~,\\
A(H,L)&=[\underline{q}(H,L),q(H,L)]~,\quad\quad &A(H)&=[\underline{q}(H),q(H)]~.
\end{alignat}
\end{proposition}
The questions of existence and uniqueness of rendezvous or average numbers, are two naturally
posed problems, which were investigated in \cite{FRpotenc} in the potential theoretic framework
on locally compact spaces. In fact, non-emptyness of the rendezvous, respectively the average
intervals means that in the above formulation \eqref{randneqchebyn} the formal lower endpoints
of the intervals do not exceed the upper endpoints (we use the convention $[a,b]=\emptyset$ if
$a>b$). While uniqueness is the same as that the respective interval reduces to one point. We
recall the following two results from \cite{FRpotenc}.

\begin{theorem}\label{th:existence} Let $X$ be a locally compact Hausdorff
space, $\emptyset \ne H \subset L \subset X$ be arbitrary, and let $k$ be any nonnegative,
symmetric kernel on $X$. Then the intervals $R_n(H,L)$, $R(H,L)$ and $A(H,L)$ are nonempty.
\end{theorem}
\begin{theorem}\label{th:unique} Let $X$ be any locally compact Hausdorff
topological space, $k$ be any l.s.c., nonnegative, symmetric kernel function, and $\emptyset\neq
K\Subset X$ compact. Then $A(K)$ consists of one single point. Furthermore, if $k$ is continuous, then even $R(K)$ consists of only one point.
\end{theorem}

When the rendezvous or the average interval $R(K)$ respectively $A(K)$ consists of one point only,
this single point is denoted by $r(K)$ or $a(K)$, respectively.

Let us close this introduction with a few remarks to explain the idea of the present approach.
Investigating the polarisation constant problem, it was found in \cite{VA} that for certain
cases the Chebyshev constants of the unit spheres $S^2$ and $S^3$ appear as polarisation
constants. The arising questions led to the systematic analysis of Chebyshev constants and also
transfinite diameters and minimal energies in the general potential theoretical framework
\cite{Bela}. Meanwhile, the second author took part in working out a general approach which
might be termed as appropriate averaging over $S^n$, to estimate the polarisation constant
\cite{PR, RS}. However, it turned out that part of the results achieved through such a
potential theory flavoured approach, were already obtained by Garc\'{\i}a-V\'azquez and Villa
\cite{GVV}, who used Gross' Theorem on the existence of rendezvous numbers successfully in the
context. That suggested that perhaps there is a way to relate the two methods, or even the
underlying theories, i.e., potential theory and rendezvous numbers. Our paper stems from this
observation.

\section{Rendezvous numbers for normed spaces}\label{sec:normsp}
In the last decade many results were obtained regarding the numerical values of rendezvous
numbers of concrete spaces and sets, see, e.g., \cite{baronti/casini/papini:2000, GVV, Lin, W1,
W2}. In this context, the following terminology was introduced.

\begin{definition} Let $(X,\|\cdot\|)$ be a normed space with unit closed
ball $B_X$ and unit sphere $S_X$. Considering $S_X$ with the norm-distance, the rendezvous numbers
of the arising metric space are called the rendezvous numbers of the normed space $X$. Accordingly,
we use the script notation
\begin{equation}\label{banachdef}
\r_n(X):=R_n(S_X)~, \qquad \r(X):=R(S_X) \qquad\text{and}\quad \a (X):= A(S_X)~.
\end{equation}
\end{definition}

It is clear that for finite dimensional normed spaces the above
notion is a special case of the general notion described in
Section \ref{sec:pot}. However, for infinite dimensional normed
spaces the metric space $S_X$ will not be locally compact, as is
usually assumed in the potential theoretic setup. There are two
ways to tackle this, one being the extension of the theory to not
necessarily locally compact but {\it metric} spaces, as is done in
\cite{FRmetr}. There we assume that the topology arises from a
metric, but relax on local compactness. Conversely, in a number of
cases it is possible to consider a {\it different topology}, in
which the metric is still lower semicontinuous, while $S_X$
becomes locally compact, hence Fuglede-type general potential
theory applies. Note that in this case the topology is \emph{not}
the metric topology, which deserves some care when working with
the theory. In particular, the average sets $A(H,L)$, referring to
regular Borel measures of the space, may be different for
different topologies. On the other hand, for a fixed kernel
$R(H,L)=[M(H,L),\overline{M}(H,L)]$ is independent of any
topology.

\begin{proposition}\label{prop:always} Let $X$ be any abstract
set, $k\geq 0$ be a symmetric function from $X\times X$ to ${\RR}_{+}\cup\{+\infty\}$, and
$\emptyset\neq H\subset L\subset X$ be arbitrary subsets. Then $R(H,L)\ne \emptyset$.
\end{proposition}
\begin{proof} The definition of rendezvous intervals, as well as the corresponding statement in
Proposition \ref{prop:reni}, are independent of the topology of the underlying space. Therefore, we
can just take the discrete topology of $X$, and note that the kernel $k$ becomes continuous, hence
l.s.c., in this topology. Thus Theorem \ref{th:existence} applies and $R(H,L)\ne\emptyset$.
\end{proof}

In view of Proposition \ref{prop:always}, we trivially obtain the following.

\begin{corollary}\label{cor:nonempty} Let $X$ be any normed space. Then the
rendezvous set $\r(X)$ is non-empty.
\end{corollary}

This seemingly contradicts to some assertions on nonexistence in the literature: the reason is
that we considered also the closure in the definition of the rendezvous and average intervals.
The above proposition shows that in the context taking the closure is also helpful. Note that
Baronti, Casini and Papini \cite{baronti/casini/papini:2000} have already considered the closed
version of the rendezvous sets, at least in normed spaces but they focused on the question of
``attaining'' the rendezvous numbers. In such investigations the geometry of Banach spaces
plays an important role.

On the other hand, uniqueness, so nicely obtained for locally compact spaces, continuous kernels
and compact sets, can not be concluded as already shown by a couple of examples in the literature
(see \cite{hinrichs}, \cite{hinrichs/wenzel}, \cite{Lin}).

\section{Average numbers for normed spaces}\label{sec:avnormsp}
For compact sets $K$ and continuous kernels it is already known
that $A(K)=R(K)$, and that there are counterexamples showing that
in general compactness is needed (compare \cite[\S 6]{FRpotenc}).
Nevertheless, the assertion remains valid in normed spaces, too.

\begin{theorem}\label{th:ARnormed} Let $X$ be any normed space.
Then we have $\a(X)=\r(X)\ne \emptyset$.
\end{theorem}

For not necessarily finite dimensional Banach spaces, we do not have the means to restrict
considerations to compactly supported measures only. Instead, we prove the following result, whose
easy consequence is the above theorem.

\begin{theorem}\label{th:ARmetricequi} Let $(Y,d)$ be a metric space. Assume
  that the kernel $k$ is positive, symmetric and bounded and that
$\{k(\cdot,y): y\in Y\}$ is uniformly equicontinuous on $(Y,d)$. Then we have $A(Y)=R(Y)\ne \emptyset$.
\end{theorem}

\begin{lemma}\label{lem:measureapprox}
Assume
  that the kernel $k$ is positive, symmetric and bounded and that
$\{k(\cdot,y): y\in Y\}$ is uniformly equicontinuous on $Y$.  Let $\mu\in
  \fMe(Y)$ and $\ve>0$ be given arbitrarily. Then there exist $m\in \NN$ and points $x_j\in Y$ ($j=1,\dots,m$)
such that the potential $U^{\nu}$ of the measure $\nu:=\frac 1m \sum_{j=1}^m \delta_{x_j}$
approximates $U^{\mu}$ within $\ve$ uniformly on $Y$.
\end{lemma}
To prove this lemma we need the following elementary result.
\begin{lemma}\label{lem:raci}
For any $\varepsilon>0$ and any finitely supported probability measure $\nu$, there exists a
probability measure of the form $\mu=\frac{1}{m}\sum_{i=1}^m \delta_{z_i}$ having the same support
as  $\nu$ and satisfying $(1-\varepsilon)\nu\leq \mu\leq (1+\varepsilon)\nu$.
\end{lemma}
\begin{proof}[Proof of Lemma \ref{lem:measureapprox}]
Without loss of generality we can assume that $k\leq 1$, and hence $U^\sigma\leq 1$ for all
probability measures $\sigma$. By the assumptions we find an $r>0$, such that $|k(x',y)-k(x'',y)| <
\ve/2$, if $d(x',x'')<r$. As $\mu$ is a \emph{regular} Borel measure, for any given $\eta
>0$ there exists $K\Subset \supp \mu \subset Y$ with $\mu(K)\geq
1-\eta$. Take now $\nu_K:=\mu_K/\|\mu_K\|$ (with $\mu_K$ being the trace of $\mu$ on $K$) so that
$\nu_K\in \fMe(K)$. Note that
$$
|U^{\mu}(x)-U^{\mu_K}(x)|\le \sup_{y\in Y} k(x,y)\cdot \|\mu-\mu_K\|  \leq \eta~,
$$
and
$$
|U^{\nu_K}(x)-U^{\mu_K}(x)|\le |U^{\nu_K}(x)|\cdot (1-\|\mu_K\|) \leq  \eta~.
$$
Consider now the covering of the compact set $K\Subset Y$ by open balls $B(y,r)$, with $y\in
K$. By compactness there exist a finite sub-covering, i.e., there exist $n\in \NN$, $y_j\in K$,
$B_j:=B(y_j,r)$ ($j=1,\dots,n$) satisfying $K\subset \bigcup_{j=1}^n B_j$. Put $D_1:=B_1$, and
for $j=2,\dots,n$ put $D_j:=B_j\setminus \bigcup_{i=1}^{j-1} B_i$, $\alpha_j:=\nu_K(D_j)\geq
0$. Clearly $\sum_{j=1}^n \alpha_j = \nu_K(K)=1$. Consider the finitely supported measure
$\sigma:=\sum_{j=1}^n \alpha_j \delta_{y_j}\in \fMe(K)$. Then we have
\begin{align*}
|U^{\nu_K}(x)-U^{\sigma}(x)|&= \Bigl|\sum_{j=1}^n \int_{D_j} k(x,y)-k(x,y_j) \dd\nu_K(y) \Bigr|
\leq\\&\leq \sum_{j=1}^n \alpha_j \sup_{y\in D_j\subset B_j} |k(x,y)-k(x,y_j)| \leq \ve/2\,.
\end{align*}
 Finally, the application of Lemma \ref{lem:raci} to $\sigma$ yields an approximating
measure $\nu:=\frac 1m \sum_{j=1}^m \delta_{x_j}$ with $m\in \NN$ and points $x_j\in K$
($j=1,\dots,m$) so that $(1-\eta)\sigma \leq \nu \leq (1+\eta) \sigma$, and thus
$$
|U^{\sigma}(x)-U^{\nu}(x)|\leq \eta |U^{\sigma}(x)| \leq\eta~.
$$
Collecting all the above, we find
$$
|U^{\mu}(x)-U^{\nu}(x)|\leq 3 \eta + \ve/2  <\ve \,,
$$
if $\eta <\ve/6$.
\end{proof}

\begin{proof}[Proof of Theorem \ref{th:ARmetricequi}] It is obvious that $A(Y)\subset R(Y)$, hence it suffices to show the converse inclusion. Let
$\mu\in \fMe(Y)$ be arbitrary and consider $A(\mu,Y)=[a,b]$, (where $a=\underline{Q}(\mu,Y)$,
$b= {Q}(\mu,Y)$). Let us take $\ve:=1/n$ and look at the measure $\nu:=\nu_n$ provided by Lemma
\ref{lem:measureapprox} to $\mu$ and $\ve$. Since the potential functions are uniformly close
to each other on $Y$, their infima and suprema are also within $\ve$: that is,
$|\underline{Q}(\mu,Y)-\underline{Q}(\nu,Y)|\leq \ve$, $|{Q}(\mu,Y)-{Q}(\nu,Y)|\leq \ve$. In
other words, $A(\nu_n,Y)\subset [a-\frac 1n,b+\frac 1n]$ and thus $\bigcap_{n=1}^{\infty}
A(\nu_n,Y) \subset [a,b] = A(\mu,Y)$. It follows that $R(Y)=\bigcap_{\nu\in\fMe^{\#}(Y)}
A(\nu,Y) \subset A(\mu,Y)$ for all $\mu\in\fMe(Y)$, hence $R(Y)\subset \bigcap_{\mu\in\fMe(Y)}
A(\mu,Y)=A(Y)$ and the theorem is proved.
\end{proof}

\begin{remark}\label{rem:equikern}
The assumptions of Theorem \ref{th:ARmetricequi} are fulfilled, for instance, when $(X,d)$ is a
metrisable topological vector space and $k(x,y)=f(d(x,y))$, where $f$ is a continuous function
(see Section \ref{sec:concrete} below).
\end{remark}

\begin{lemma}\label{lem:convex} Let $X$ be a (not necessarily
locally compact) Hausdorff topological vector space, $k$ a l.s.c., nonnegative, symmetric and
convex kernel function on $X\times X$, and $\mu\in \fMe(X)$. Then the potential function $U^{\mu}$
is convex.
\end{lemma}
\begin{proof} Take any $x,y,z\in X$ with $z=\alpha x +(1-\alpha)
y$, where $0\leq \alpha \leq 1$. We then have
\begin{align*}
U^{\mu}(z)&=\int_X k(z,w) \dd\mu(w) \leq \int_X (\alpha k(x,w) +(1-\alpha)k(y,w))
\dd\mu(w)=\\
&=\alpha U^{\mu}(x)+(1-\alpha) U^{\mu}(y)~,
\end{align*}
and that was to be proved.
\end{proof}

\begin{lemma}\label{lem:boundary} Let $X$ be a (not necessarily
locally compact) Hausdorff topological vector space, $k$ a l.s.c., nonnegative, symmetric and
convex kernel function on $X\times X$, $H\subset X$ a bounded set, and $\partial H$ be its
boundary. Then for any $\mu\in \fMe(X)$ the potential function $U^{\mu}$ satisfies $\sup_H U^{\mu}
= \sup_{\partial H} U^{\mu}$.
\end{lemma}
\begin{proof} We are to show $\sup_H U^{\mu} \leq \sup_{\partial H} U^{\mu}$,
the other direction being obvious. Let now $x\in H$ be arbitrary: we show that $U^{\mu}(x) \leq
\sup_{\partial H} U^{\mu}$. Draw any straight line $\ell$ through $x$. Since $H$ is bounded, both
closed half-lines of $\ell$, starting from $x$, contain boundary points of $H$; that is, if these
points are $y, z \in \ell$, then $x\in [y,z]$ with $y,z \in \partial H$. According to Lemma
\ref{lem:convex}, the potential is convex, which immediately yields $U^{\mu}(x)\leq \max
(U^{\mu}(y),U^{\mu}(z))\leq \sup_{\partial H} U^{\mu}$.
\end{proof}

\begin{corollary} If $X$ is a normed space, then
$q(S_X,S_X)=q(S_X,B_X)$ and $\overline{M}(S_X,S_X)=\overline{M}(S_X,B_X)$.
\end{corollary}

\begin{remark} Note that for the other endpoints of the average intervals generally we may have strict
inequality: $\underline{q}(S_X,B_X)<\underline{q}(S_X,S_X)$. As
will be seen in Theorem \ref{th:lpball}, for any $1< p <+\infty$
$\r(\ell_p)=2^{1/p}>1$. However, for any measure
$\mu\in\fMe(S_X)$, it is clear that ${\mathbf 0} \in B_X$ satisfies
$U^{\mu}({\mathbf 0})=\int_{S_X} 1 \dd\mu =1$, hence
$\underline{q}(S_X,B_X)\leq 1$. In fact, $\underline{q}(S_X,B_X)=
1$ is also true, since for any set $H$ and for two points $x,y\in
H$, the corresponding measure $\nu:=\frac 12(\delta_x+\delta_y)$
always provides, by the triangle inequality,
$\underline{Q}(\nu,B_X) \geq \|x-y\|/2$, which can be as large as
half the diameter.
\end{remark}

\begin{remark}\label{rem:diam2diam}
It is straightforward to show that $[M(H),\overline{M}(H)]=R(H)\subset
[\frac{1}{2}\diam(H),\diam(H)]$ (see, e.g., \cite{Gross}). Indeed, $\overline{M}(H)\leq \diam
(H)$ is trivial, while the lower estimate $\frac{1}{2}\diam(H)\leq M(H)$ is essentially
contained in the previous remark.
\end{remark}
\section{Rendezvous numbers for $L_{p}$ spaces}\label{sec:concrete}
We already know that the rendezvous interval of a Banach space is not empty. Here we identify the
rendezvous, hence the average intervals of the $L_p$ spaces.

Let $(\Omega,\mathcal{M},\mu)$ be a measure space. To complement the whole scale $1\leq p<+\infty$, we
consider $L_p:=L_p(\Omega,\mathcal{M},\mu)$ when $0<p<1$ as well. In this case $L_p$ will not be a
Banach space, but if endowed with the metric
\begin{equation}\label{eq:metrdef}
d(f,g)=\int_\Omega |f-g|^p\dd\mu~,
\end{equation}
 it is a complete, metrisable, topological
vector space (of course we have to identify functions coinciding on $\mu$-null sets). First we
calculate the rendezvous number with respect to the symmetric function $\|f-g\|:=d(f,g)^{1/p}$
instead of the metric $d$, this fits well together with the case $p\geq 1$. Of course, now
$S_{L_p}$ denotes the ``unit sphere'' with respect to $\|\cdot\|$ for all $0<p<+\infty$.

\begin{theorem}\label{th:lpball}
Let $0< p<+\infty$ be arbitrary and consider
$L_p(\Omega,\mathcal{M},\mu)$ over either $\RR$ or $\CC$. If $L_p$
is infinite dimensional, we have $a(S_{L_p})=r(S_{L_p})=2^{1/p}$.
\end{theorem}
\begin{proof}
The following applies for both the complex- and real valued cases, hence we do
not mention the underlying number field any more. Further we write briefly
$L_p$ instead of $L_p(\Omega,\mathcal{M},\mu)$.  We already know that
$R(S_{L_p})$ is nonempty (see Proposition \ref{prop:always}) and is a compact
interval (indeed, in case $p\geq 1$ it is a subset of the interval $[1,2]$,
see Remark \ref{rem:diam2diam}). Moreover, we have
$R(S_{L_p})=[M(S_{L_p}),\overline{M}(S_{L_p})]$, and by Theorem
\ref{th:ARnormed}, Theorem \ref{th:ARmetricequi} and Remark \ref{rem:equikern} we have
$A(S_{L_p})=R(S_{L_p})$. Therefore we need only to show that $M(S_{L_p}) \geq
2^{1/p}$ and that $\overline{M}(S_{L_p}) \leq 2^{1/p}$.

Since $L_p$ is infinite dimensional, there exist $w_j\in S_{L_p}$,
$j\in \NN$ such that the sets $A_j:=\{x:\:x\in \Omega,\:w_j(x)\neq
0 \}$ are pairwise disjoint. For any function $g\in L_p$ let us
introduce the notation $\|g\|_j:=\|\chi_{A_j}g\|$.

\emph{Part 1: $M(S_{L_p}) \geq 2^{1/p}$.} With the functions $w_j$ and for any $f\in S_{L_p}$ we have
\begin{align}\label{Mestimate1}
\sum_{j=1}^n \|f-w_j\|&=\sum_{j=1}^n \Bigl( \sum_{k=1}^\infty\|f-w_{j}\|^p_k \Bigr)^{\frac
1p}=\sum_{j=1}^n \Bigl( \|f-w_j\|^p_j +
\sum_{k=1, k\ne j}^\infty \|f\|^p_k\Bigr)^{\frac 1p}=\notag \\
&= \sum_{j=1}^n \left( \|w_j-f\|^p_j+1- \|f\|^p_j \right)^{\frac 1p}~.
\end{align}
Now we distinguish between the cases $p<1$ and $p\geq 1$. First let $p<1$, then using $\|f\|_j\leq
1$ and $f\in S_{L_p}$, we can continue \eqref{Mestimate1}
\begin{align*}
&\sum_{j=1}^n \|f-w_j\|= \sum_{j=1}^n \left( 1- \|f\|^p_j + \|w_j-f\|^p_j \right)^{\frac 1p}\geq\\
&\quad\geq\sum_{j=1}^n \left( 1-  \|f\|^p_j+\left|\|w_j\|_j^p-\|f\|^p_j\right|\right)^{\frac 1p}=\sum_{j=1}^n \left( 1-  \|f\|^p_j+\left(1-\|f\|^p_j\right)\right)^{\frac 1p}=\\
&\quad=2^{1/p}\sum_{j=1}^n \left( 1-  \|f\|^p_j\right)^{\frac 1p}\geq 2^{1/p}n
\frac{(n-1)^{\frac 1p}}{n^{\frac 1p}}~,
\end{align*}
using again $f\in S_{L_p}$ and the convexity of the function $x\mapsto x^{1/p}$ in the last step.
Second, let $p\geq 1$, then we can write
\begin{align}\label{Mestimate2}\notag
\sum_{j=1}^n \|f-w_j\|&= \sum_{j=1}^n \left( 1- \|f\|^p_j + \|w_j-f\|^p_j \right)^{\frac
1p}\geq\\
&\geq\sum_{j=1}^n \left( 1-  \|f\|^p_j+\left(1-\|f\|_j\right)^{p}\right)^{\frac 1p}=\notag\\
&= \sum_{j=1}^n \left( 1-  \|f\|^p_j-\left(1-\|f\|_j\right)^{p}+2\left(1-\|f\|^p_j\right)^{p}\right)^{\frac 1p}\geq\\
&\geq\sum_{j=1}^n \left(2\Bigl(1-\|f\|_j\Bigr)^{p}\right)^{\frac 1p}=2^{\frac 1p}\sum_{j=1}^n
\left(1-\|f\|_j\right)\geq 2^{\frac 1p}( n - n^{\frac 1q})~,\notag
\end{align}
where $\frac 1p +\frac 1q$ and using again $f\in S_{L_p}$ and H\"older's inequality in the last
step. We see that in both cases for the given $n$-point distribution (concentrated on the
$w_j$s, $j=1,\dots,n$) the corresponding potential is minorised by the right hand sides divided
by $n$, hence we find $M_n(S_{L_p}) \geq 2^{1/p}+o(1)$ (as $n\to \infty$) and $M(S_{L_p})\geq
2^{1/p}$ follows. (Note that, e.g., \cite{GVV} calculates even the exact value of
$\r(\ell_1^n(\CC))$, but here the task is a little bit different.)

\emph{Part 2: $\overline{M}(S_{L_p}) \leq 2^{1/p}$.} Let $n\in \NN$ and consider the same functions
$w_j$, $j\in \NN$ as in the first part. Then for any point $f \in S_{L_p}$ and for any given
parameter $\eta>0$ we have, by \eqref{Mestimate1}, for $p\geq 1$
\begin{align}\label{Mestimateup}
\sum_{j=1}^n \|f-w_j\| &\leq \sum_{j=1}^n \left( 1 + (1+\|f\|_j)^p \right)^{\frac 1p} \leq\\
&\leq \sum_{j~: \|f\|_j>\eta} \left( 1+ 2^p \right)^{\frac 1p} + \sum_{j~: \|f\|_j\leq \eta}
\left( 1+ (1+\eta)^p \right)^{\frac 1p}\leq\notag\\
&\leq \frac{1}{\eta^p} \left(1+2^p\right)^{\frac 1p} + n 2^{\frac 1p} (1+\eta)~,\notag
\end{align}
and for $p<1$
\begin{align}\label{Mestimateup2}
\sum_{j=1}^n \|f-w_j\|& \leq \sum_{j=1}^n \left( 1 + (1+\|f\|^p_j) \right)^{\frac 1p} \leq\notag\\
&\leq \sum_{j: \|f\|_j>\eta} 3^{\frac 1p} + \sum_{j: \|f\|_j\leq \eta} \left( 1+ (1+\eta^p)
\right)^{\frac 1p}\leq \frac{3^{\frac 1p}}{\eta^p}  + n 2^{\frac 1p} (1+\eta^p)^{\frac{1}{p}}~.
\end{align}
 Choosing, e.g., $\eta:=n^{-1/{2p}}$, we obtain the estimate
 $U^{\nu}(f)\leq
2^{1/p} + \mathrm{o}_p(1)$ (as $n\rightarrow+\infty$, $\forall
f\in S_{L_p}$) for the measure $\nu:=\frac 1n\sum_{j=1}^n
\delta_{w_j}$. It follows that for the given $n$-point
distribution $\nu$ we have $\overline{M}_n(S_{L_p})\leq
Q(\nu,S_{L_p})=2^{1/p} + \mathrm{o}_p(1)\to 2^{1/p}$ ($n\to
\infty$), and so even $\overline{M}(S_{L_p}) \leq 2^{1/p}$.
\end{proof}

Note that for $p\geq 1$ already Lin \cite{Lin} showed that
$\r(\ell_p)\subseteq\{2^{1/p}\}$ for ``strict'' rendezvous numbers
(actually by a similar argument). So this and the non-emptyness of
the rendezvous interval (Corollary \ref{cor:nonempty}) give the
above result for $\ell_p$ ($1\leq p< +\infty$).
\begin{corollary}[\bf Wolf, Lin\rm] Let ${\mathcal H}$  be an infinite dimensional Hilbert space over any of the number fields $\RR$ or $\CC$. Then
we have $\a({\mathcal H})=\r({\mathcal H})=\{\sqrt{2}\}$.
\end{corollary}

\begin{remark} In the above proof we actually used the same point
distribution in both parts of the proof, therefore we have proved
the existence of $\ve$-\emph{quasi-invariant measures}. We say
that there exist \emph{$\varepsilon$-quasi-invariant measures}
cf.~\cite[Definition 5.9]{FRmetr} for the kernel $k$ on $S$, if
for all $\ve>0$ there is $\nu\in\fMe(S)$ satisfying
$Q(\nu,S)-\underline{Q}(\nu,S) \leq \ve$. Generally, if $S$ is
compact, by  \ws-compactness we obtain the existence of a true
invariant measure $\mu\in\fMe(S)$, i.e., whose potential $U^{\mu}$
is constant on $S$. Of course, for $S=S_{\ell_p}$ this is not the
case. But as we saw above there exist
$\varepsilon$-quasi-invariant measures, and by \cite[Proposition
5.11]{FRmetr} it is already enough to conclude that the average
interval reduces to one point.
\end{remark}

\begin{theorem}\label{th:lprealball}
Let $0<p<1$ be arbitrary and  $L_p(\Omega,\mathcal{M},\mu)$ be the vector space of $p$-integrable functions endowed with the metric
$d$ defined in \eqref{eq:metrdef}. For the rendezvous interval of the unit ball $S_{L_p}$ of
$L_p$ we have $A(S_{L_p})=R(S_{L_p})=2$.
\end{theorem}
\begin{proof}
By Remark \ref{rem:diam2diam} we know $\overline{M}(S_{L_p})\leq 2$. So we only need to estimate
the lower endpoint of $R(S_{L_p})$ from below. This can be done analogously to
\eqref{Mestimate1} and \eqref{Mestimate2} by considering the functions $w_j$ and sets $A_j$ used in the proof of Theorem \ref{th:lpball}. Let us further use
the abbreviation $\nm f\nm_j:=\int_{A_j} |f|^p\dd\mu$. For an arbitrary $f\in S_{L_p}$ we can write
\begin{align*}
\sum_{j=1}^n d(f,w_j)&=\sum_{j=1}^n \left( \sum_{k=1}^\infty \nm f-w_{j}\nm_k\right)=\sum_{j=1}^n
\left( \nm f-w_j\nm_j +
\sum_{k=1, k\ne j}^\infty \nm f\nm_k \right)=\\
&= \sum_{j=1}^n \left( 1- \nm f\nm_j + \nm w_j-f\nm_j \right)\geq\\
&\geq\sum_{j=1}^n \left( 1- \nm f\nm_j + 1-\nm f\nm_j \right)=2\sum_{j=1}^n
\left( 1-\nm f\nm_j \right)\geq2
(n-1)~.
\end{align*}
Therefore $M_n(S_{L_p}) \geq 2-1/n$, hence $M(S_{L_p})\geq 2$.
\end{proof}

It was already pointed out by Lin in \cite{Lin} that $\lim_{n\rightarrow\infty}
r(S_{\ell_p^n})=2^{1/p}$, if $1\leq p<+\infty$. The following result, inspired by an analogous
argument of Garc\'{\i}a-V\'azquez and Villa \cite{GVV}, explains this phenomenon in view of Theorem
\ref{th:lpball} above.

\begin{theorem}\label{thm:xnwithnorm}
Let $X$ be a normed space and $X_n$ an increasing sequence of subspaces such that
$\bigcup_{n=1}^\infty X_n$ is dense in $X$. Suppose that $(X_n,\|\cdot\|_n)$ is a normed space and
$\rho_n\in R(S_{(X_n,\|\cdot\|_n)})$. Assume that
$$
\lim_{n\rightarrow \infty} \sup_{x\in X_n\cap S_X} |1-\|x\|_n|=0~.
$$
Then any accumulation point $\rho$ of the sequence $\rho_n$ belongs to $R(S_X)$.
\end{theorem}
\begin{proof}
Let $\rho$ be an accumulation point of the sequence $\rho_n$.  Assume without loss of
generality that $\rho_n\to \rho$. Let $\varepsilon>0$ be given. Then for sufficiently large
$n\geq n_0(\varepsilon)$ we have $[\rho_n-\varepsilon,\rho_n+\varepsilon]\subset
[\rho-2\varepsilon,\rho+2\varepsilon]$.

By definition $\rho_n\in \r_m(X_n)$ for all $m\in \NN$. Let $m\in \NN$ be fixed, and
$x_1,\dots, x_m\in S_X$ be arbitrary. Take any $y_j\in X_n\cap S_X$ with
$\|y_j-x_j\|\leq\varepsilon$. Such $y_j$ exists in view of the denseness of
$\bigcup_{n=1}^\infty X_n$ in $X$. By assumption, we have $|\|z\|-\|z\|_n|\leq \varepsilon
\|z\|$ for all sufficiently large $n\geq n_1\geq n_0$ and all $z\in X_n$. In particular,
$|1-\|y_j\|_n|\leq \varepsilon$. By definition of the rendezvous interval $\r_m(X_n)$, we find
$z_n\in S_{(X_n,\|\cdot\|_n)}$ satisfying
\begin{equation}\label{eq:randibol}
\frac{1}{m}\sum_{j=1}^m
\Bigl\|\frac{y_j}{\|y_j\|_n}-z_n\Bigr\|_n\in[\rho_n-\varepsilon,\rho_n+\varepsilon]\subset
[\rho-2\varepsilon,\rho+2\varepsilon]~.
\end{equation}
According to the above, for all $n\geq n_1$ we have
\begin{align*}
\left\|z_n-\frac{z_n}{\|z_n\|}\right\|_n&=\left|1-\frac{1}{\|z_n\|}\right|\cdot\|z_n\|_n =
\left|1-\frac{\|z_n\|_n}{\|z_n\|}\right| \leq \varepsilon~, \intertext{and}
\left\|y_j-\frac{y_j}{\|y_j\|_n}\right\|_n&=\left|1-\frac{1}{\|y_j\|_n}\right| \cdot \|y_j\|_n =
\left|{\|y_j\|_n}-1\right|\leq \varepsilon~.
\end{align*}
Using these two inequalities in \eqref{eq:randibol} we obtain
$$
\frac{1}{m}\sum_{j=1}^m \Bigl\|y_j-\frac{z_n}{\|z_n\|}\Bigr\|_n\in
[\rho-4\varepsilon,\rho+4\varepsilon]~.
$$
For $n>n_1$ we also know
$$
\left|\Bigl\|y_j-\frac{z_n}{\|z_n\|}\Bigr\|-\Bigl\|y_j-\frac{z_n}{\|z_n\|}\Bigr\|_n\right|\leq
\varepsilon \cdot \Bigl\|y_j-\frac{z_n}{\|z_n\|}\Bigr\|\leq 2\varepsilon~,
$$
therefore we can write
$$
\frac{1}{m}\sum_{j=1}^m \Bigl\|y_j-\frac{z_n}{\|z_n\|}\Bigr\|\in
[\rho-6\varepsilon,\rho+6\varepsilon]\quad\mbox{and}\quad\frac{1}{m}\sum_{j=1}^m
\Bigl\|x_j-\frac{z_n}{\|z_n\|}\Bigr\|\in [\rho-7\varepsilon,\rho+7\varepsilon]~.
$$
This shows $\rho\in R_m(S_X)$, which being valid for all $m$, gives $\rho\in R(S_X)$.
\end{proof}

This theorem immediately gives the following corollary.

\begin{corollary}\label{cor:xn}
Let $X$ be a normed space and $X_n$ an increasing sequence of finite dimensional subspaces such
that $\bigcup_{n=1}^\infty X_n$ is dense in $X$. Let $\rho$ be an accumulation point of the
sequence $r(S_{X_n})$ ($r(S_{X_n})$ exists uniquely by the compactness of $S_{X_n}$). Then $\rho\in
\r(X)$.
\end{corollary}

\section{Chebyshev centres, entropy and rendezvous numbers} \label{ss:center}

\begin{definition} Let $K\subset X$ be a compact, convex subset of
some normed space $X$, with $d$ being the metric induced by the norm. Then the Chebyshev centre
$c:=c(K)\in K$ and the Chebyshev out-radius $\rho:=\rho(K)$ of the set $K$ are the centre and the
radius, respectively, of the closed ball $\overline{B}:=\overline{B}(c,\rho)$ of minimal radius
with $c\in K\subset\overline{B}(c,\rho)$.
\end{definition}

Clearly, for any compact $K$ such a minimal ball always exists, and for convex sets it is even
unique. Note that it is important in the definition that $c$ be chosen within $K$; for
discussion see \cite{BK} and \cite{CMY, G}.

Quoting private communication from Esther and George Szekeres, in \cite{CMY} Cleary, Morris and
Yost present the following beautiful result with a nice, direct elementary proof. Here we
present our even shorter version relying on the potential theoretical background developed.

\begin{theorem}[\bf E.~and G.~Szekeres\rm] Let $K\subset X$ be a compact,
convex subset of some normed space $X$, with $d$ being the metric induced by the norm. Then the
Chebyshev radius and the rendezvous number of the set $K$ equal: $r(K)=\rho(K)$.
\end{theorem}

\begin{proof} Existence and uniqueness of
$R(K)=\{r(K)\}$ and also the equality $R(K)=A(K)$ are already known from Theorems
\ref{th:existence} and \ref{th:unique}. Further, if $c$ is a Chebyshev centre of $K$, then to all
points $x\in K$ we have $\|x-c\|\leq \rho$. Hence for any probability measure $\mu\in\fMe(K)$ also
the potential satisfies $U^{\mu}(c)\le \rho(K)$. It follows that $\underline{Q}(\mu,K)\leq \rho$
for all $\mu\in\fMe(K)$, hence also $ r(K)=\underline{q}(K)\leq \rho(K)$.

Conversely, for $\ve>0$ let $y_j\in K$ be arbitrary points ($j=1,\dots,n$) satisfying $Q(\nu,K)
< \overline{M}(K)+\ve$ with $\nu:=\frac 1n \sum_{j=1}^n \delta_{y_j}$.  As $K$ is convex, it
contains the convex combination $y:=\frac 1n \sum_{j=1}^n y_j \in K$ of the given points. Thus
by the convexity of the norm for arbitrary $x\in K$ the estimate $\|y-x\|\leq \frac 1n
\sum_{j=1}^n \|y_j-x\|=U^{\nu}(x) \leq \overline{M}(K)+\ve $ holds. Hence
$\overline{B}(y,\overline{M}(K)+\ve)$ covers $K$ and $\rho(K)\leq \overline{M}(K)+\ve$, which
implies also $\rho(K)\le \overline{M}(K)=q(K)=r(K)$.
\end{proof}

Recall that for a positive number $t>0$ and a set $H\subset X$ of a metric space $X$ with metric
$d$ the $t$-covering number $N(t,H)$ is defined as
$$
N(t,H):=\min \{n\in \NN ~:~ \exists y_j\in H~ (j=1,\dots,n) \:\:\text{such that}~ H \subset
\bigcup_{j=1}^n B(y_j,t)\}~.
$$
If there is no finite set of balls of radius $t$ which can cover the set $H$, then we say that
$N(t,H)=+\infty$. Similarly, if $H,L\subset X$, then
$$
N(t,H,L):=\min \{n\in \NN ~:~ \exists y_j\in H~ (j=1,\dots,n) \:\:\text{such that}~ L \subset
\bigcup_{j=1}^n B(y_j,t)\}
$$
with $\min \emptyset =+\infty$ being in effect again. The next observation is almost obvious.
\begin{proposition}\label{pr:entropy} Let $t>0$ and $H,L
\subset X$. We have $ \overline{M}_n(H,L) \geq t$ for all $n < N(t,H,L)$. In particular, if
$N(t,H,L)=+\infty$, then $t \le \overline{M}(H,L)=\sup R(H,L)$.
\end{proposition}
\begin{proof} Since $n<N(t,H,L)$, for any system of points
$y_j\in H ~ (j=1,\dots,n) $, there exists some point $x\in L$ so that $d(x,y_j)\geq t$ for all
$j=1,\dots,n$. Therefore, $\sup_{x\in L} \sum_{j=1}^n d(x,y_j) \geq nt$ holds
  for all systems
of $n$ points, whence $\overline{M}_n(H,L)\geq t$, and \eqref{randneqchebyn} concludes the proof.
\end{proof}

Recall that a set $H\subset X$ is called \emph{totally bounded}, if $N(t,H)<+\infty$ for all
$t>0$. In Banach spaces this is the same as the conditional compactness of $H$, i.e., that
$\overline{H}\Subset X$ is compact set. The proposition shows that for subsets $H$ which are
not totally bounded, there is always a positive lower bound of $\overline{M}(H,H)$. The above
proposition, however easy, provides an essential help in describing some rendezvous numbers.
For instance, there is an elegant interpretation of the following result.

\begin{theorem}
\label{th:CK} Let $K$ be a compact Hausdorff topological space without isolated points, and
consider $C(K)$ the Banach space of real-  or complex-valued continuous functions over $K$. Then we
have $\overline{M}(S_{C(K)})=2$.
\end{theorem}
\begin{proof}
Denote by $S$ the unit sphere of $C(K)$.  We show that $N(t,S,S)=+\infty$ for $0<t<2$, then by
Proposition \ref{pr:entropy}  $\overline{M}(S)\geq t$ hence $\overline{M}(S)\geq 2$ will
follow. Then, by Remark \ref{rem:diam2diam}, we must actually have $\overline{M}(S)=2$.

So let $0<t<2$ and $f_1,\dots f_m\in S$. Further let $x_j\in K$ be one of the maximum points of
$|f_j|$, i.e., $|f_j(x_j)|=1$ and take $\varepsilon>0$ small such that $t+\varepsilon<2$. By
continuity, there exist neighbourhoods $G_j$ of $x_j$ with $|f_j(x_j)-f_j(y)|<\varepsilon$,
$(\forall y\in G_j)$, for all $j=1,\dots, m$. Take $y_j\in G_j$ distinct points ($x_j$ is not
an isolated point!). By Tietze's Theorem  there exists a continuous function $g\in S$, such
that $g(y_j)=-f_j(x_j)$, for all $j=1,\dots, m$. But then $S$ can not be covered by the balls
$B(f_j,t)$, because this particular $g$ is not covered. Thus we conclude $N(t,S,S)=+\infty$.
\end{proof}
The above result is already present in Garc\'{\i}a-V\'azquez,
Villa \cite{GVV} and Lin \cite{Lin}, where the authors determine
the rendezvous interval of $C(K)$. Their proofs follow the same
line, we included it for the sake of illustration of the role of
entropy.

The real-valued case in the following theorem is due to Wolf \cite{W2}, see also Lin \cite{Lin}.

\begin{theorem}Let $c_0$ denote the Banach space of real or complex valued
null-sequences. Then we have
\begin{equation}\label{corandi}
\a(c_0)=\r(c_0)=\bigl[1, \sigma]~,
\end{equation}
where $\sigma=3/2$ or $\sigma=\frac 13 + \frac{2\sqrt{3}}{\pi}$ in the case of $\RR$ respectively
$\CC$-valued sequences.
\end{theorem}
\begin{proof}
In the real case Wolf showed that $r(S_{\ell_\infty^n})=\sigma$ (\cite[Proposition 1]{W2}),
while the corresponding equality in the complex case is due to Garc\'{\i}a-V\'azquez and Villa.
So by Corollary \ref{cor:xn} we see that $\sigma\in \r(c_0)$. Applying the same idea as in
\cite[Theorem 5]{GVV}, we can split the space as $c_0=\KK\times c_0$, where $\KK$ is either the
complex or the real scalar field. Consider the measure $\mu=\lambda\otimes \delta_0$, where
$\lambda$ is the normalised Haar measure on $S_\KK$ and $\delta_0$ is the Dirac measure on
$c_0$ at the constant $0$ sequence. Clearly $\mu$ is supported in $S_{c_0}$. In case of
real-valued sequences, Wolf essentially showed $Q(\mu,S_{c_0})\leq \sigma$ (see proof of
Proposition 1 in \cite{W2}). Moreover, one can repeat the arguments from \cite{GVV} to see that
$Q(\mu,S_{c_0})\leq \sigma$ in the complex case, too. So in both cases we have $q(S_{c_0})\leq
\sigma = \sigma(\KK)$ and as by Theorem \ref{th:ARnormed} we know
$\overline{M}(S_{c_0})=q(S_{c_0})$, we find $\overline{M}(S_{c_0})\leq \sigma$. Because
$\sigma\in\r(c_0)$, the only possibility is $\overline{M}(S_{c_0})=\sigma$.

To calculate the lower endpoint of the rendezvous interval, let now $m\in\NN$ be fixed and
$x_1,\dots,x_m\in S_{c_0}$ be arbitrary. For $\varepsilon>0$ take $n_0\in\NN$ be so large that
$|x_j(n)|<\varepsilon$ whenever $n\geq n_0$ for all $j=1,\dots,m$. Now let $z$ be the element of
$c_0$ being almost completely $0$ but $1$ at the $n_0$th coordinate. Then
$$
\frac{1}{m}\sum_{j=1}^m \|x_j-z\|\leq 1+\varepsilon,
$$
so $M_m(S_{c_0})\leq 1$, and therefore $M(S_{c_0})\leq 1$. But then by Remark \ref{rem:diam2diam} we have $M(S_{c_0})=1$. We have
calculated the lower and upper endpoints of the rendezvous (and the average) intervals to arrive at
the assertion.
\end{proof}

\parindent0pt

\end{document}